\documentclass[leqno,12pt]{article}
\usepackage{amsmath}
\usepackage{amssymb}
\title{On the existence of K\"ahler metrics of constant scalar curvature}
\author{\textsc{Kenji Tsuboi}}
\date{}
\newtheorem{The}{Theorem}[section]

\newtheorem{rem}[The]{Remark}

\newtheorem{asmp}[The]{Assumption}

\def\C{{\Bbb C}}
\def\R{{\Bbb R}}
\def\Z{{\Bbb Z}}
\def\uone{\mbox{U}(1)}
\def\comb#1#2{\left(\begin{array}{c} #1\\ #2\end{array}\right)} 
\def\am#1{\mbox{{\rm Aut}}(#1)}
\def\hm#1{\frak h(#1)}
\def\cp#1{\Bbb C\Bbb P^{\,#1}}
\def\hz#1#2{H^{#1}(#2;\Bbb Z)}
\def\hr#1#2{H^{#1}(#2;\Bbb R)}
\def\GL#1{\mbox{GL}(#1,\Bbb C)}
\def\gll#1{\frak{gl}(#1,\Bbb C)}
\def\tu{\tilde{u}}
\def\tv{\tilde{v}}
\def\tw{\tilde{w}}
\begin{document}
\maketitle
\footnote[0]
   {2000\textit{Mathematics Subject Classification}.
   Primary 53C25; Secondary 53C55.}
\footnote[0]{The author is deeply grateful to the referee for valuable information.}
\begin{abstract}
For certain compact complex Fano manifolds $M$ with reductive Lie algebras of holomorphic 
vector fields, we determine the analytic subvariety of the second cohomology group of $M$ consisting 
of K\"ahler classes whose Bando-Calabi-Futaki character vanishes.  
Then a K\"ahler class contains a K\"ahler metric of constant scalar curvature if and only if 
the K\"ahler class is contained in the analytic subvariety.  
On examination of the analytic subvariety, it is shown that $M$ admits infinitely many nonhomothetic 
K\"ahler classes containing K\"ahler metrics of constant scalar curvature but does not admit any 
K\"ahler-Einstein metric.
\end{abstract}
\section{Introduction}
The question of whether a manifold admits a Riemannian metric of constant scalar curvature or not is 
a classical problem.  For any real closed manifold $M$ of dimension greater than two, Kazdan and Warner 
{\rm \cite{KW}} proved that $M$ admits at least a Riemannian metric of negative constant scalar curvature.  
On the other hand, there exists an obstruction to the existence of K\"ahler metrics of constant 
scalar curvature.  Indeed, let $M$ be an $m$-dimensional compact complex manifold.  
Denote by $\am{M}$ the complex Lie group consisting of all biholomorphic automorphisms of $M$ and 
by $\hm{M}$ its Lie algebra consisting of all 
holomorphic vector fields on $M$.  The Lie algebra $\hm{M}$ is called reductive if $\hm{M}$ is the 
complexification of the Lie algebra of a compact subgroup of $\am{M}$.  
In {\rm \cite{M}}, Matsushima proved that $\hm{M}$ is the complexification of the real Lie algebra consisting 
of all infinitesimal isometries of $M$, and hence $\hm{M}$ is reductive, if $M$ admits a 
K\"ahler-Einstein metric.  Generalizing the result of Matsushima, Lichnerowicz proved in {\rm \cite{L1}, \cite{L2}} 
that $\hm{M}$ must satisfy a certain condition if $M$ admits a K\"ahler metric of constant 
scalar curvature.  (For details see also {\rm \cite[Theorem 6.1]{K}}.)  When $M$ is a compact simply connected 
K\"ahler manifold, the condition of Lichnerowicz is equivalent to that of Matsushima.  
For example, the one point blow-up of $\cp{2}$ does not satisfy the condition (see {\rm \cite[p.100]{F3}}) 
and hence does not admit any K\"ahler metric of constant scalar curvature.  
Thus the problem to solve is whether $M$ with reductive $\hm{M}$ admits a K\"ahler metric 
of constant scalar curvature or not.

Generalizing the result of Futaki {\rm \cite{F1}},  Bando {\rm \cite{B}}, Calabi {\rm \cite{C}} and Futaki 
{\rm \cite{F2}} give an obstruction to the existence of a K\"ahler metric of constant scalar curvature 
whose K\"ahler form is contained in some particular K\"ahler class.  Let $\Omega$ be a K\"ahler class, 
$\omega\in\Omega$ a K\"ahler form and $s_{\omega}$ the scalar curvature of $\omega$.  
Let $c_1(M)\in\hz{2}{M}$ be the first Chern class of $M$ and set
$$
\mu_{\Omega}=\frac{(\Omega^{m-1}\cup c_1(M))[M]}{\Omega^m[M]}\,,
$$
where $[M]$ denotes the fundamental cycle of $M$.  Then there exists uniquely a smooth function $h_{\omega}$ up 
to constant such that
$$
s_{\omega}-m\mu_{\Omega}=\triangle_{\omega} h_{\omega}\,,
$$
and the integral
$$
f_{\Omega}(X)=\int_{M}Xh_{\omega}\omega^m
$$
is defined for $X\in\hm{M}$.  This integral $f_{\Omega}(X)$ is independent of the choice of K\"ahler 
forms $\omega\in\Omega$.  Moreover, $f_{\Omega}:\hm{M}\longrightarrow\Bbb C$ is a Lie algebra character and 
$f_{\Omega}$ vanishes if $\Omega$ contains a K\"ahler metric of constant scalar curvature.  
The character $f_{\Omega}$ is called the Bando-Calabi-Futaki character or the Futaki invariant.

When $\Omega$ is a Hodge class and a holomorphic line bundle $L$ with $c_1(L)=\Omega$ admits a lifting 
of the $\Omega$-preserving action of a subgroup $G$ of $\am{M}$, in {\rm \cite{N2}} Nakagawa gives a lifting of 
the Lie algebra character $f_{\Omega}$ to a group character $G\longrightarrow \C/(\Z + \mu_{\Omega}\Z)$ 
by using the results in {\rm \cite{T}} and {\rm \cite{FM}}.

Assume that there exists an inclusion $\iota:\uone\longrightarrow\am{M}$ and that $\Omega$ is equal to 
the first Chern class of a holomorphic $\uone$-line bundle $L$ over $M$.  
For any integer $p\geq 2$ let $Y$ denote the element $2\pi\sqrt{-1}$ of the Lie algebra of 
$\uone$ and set
\begin{equation}\label{xp}
X=\iota_{*}Y\in\hm{M}\,,\;X_p=\frac1p X\in\hm{M}\,,\;g_p=\exp X_p\in\am{M}\,.
\end{equation}
Then the order of $g_p$ is $p$.  We assume that the next assumption is satisfied.  
(See Assumption 2.2 and Lemma 2.3 in {\rm \cite{FT}}.)
\begin{asmp}\label{asp}
Assume that the fixed point set of $g_p^k$ for $1\leq k\leq p-1$ is independent of $k$ 
and that the connected components $N_1,\,\cdots,\,N_n$ of the fixed point set, which are compact 
complex submanifolds of $M$, have cell decompositions with no codimension one cells.
\end{asmp}

Let $\alpha_p$ denote the primitive $p$-th root of unity defined by
$$
\alpha_p=e^{2\pi\sqrt{-1}/p}
$$
hereafter.  Suppose that $g_p^k$ acts on $K_M^{-1}|_{N_i}$ via multiplication 
by $\alpha_p^{kr_i}$ and acts on $L|_{N_i}$ via multiplication by $\alpha_p^{k\kappa_i}$.  
Suppose moreover that the normal bundle $\nu(N_i,M)$ is decomposed into the direct sum of subbundles
$$
\nu(N_i,M)=\oplus_j\nu(N_i,\theta_j)\,,
$$
where $g_p^k$ acts on $\nu(N_i,\theta_j)$ via multiplication by $e^{\sqrt{-1}\theta_j}$.  Then a cohomology class 
$\Phi(\nu(N_i,M))$ is defined by
$$
\Phi(\nu(N_i,M))=\prod_j\prod_{k=1}^{R_j}\frac1{1-e^{-x_k-\sqrt{-1}\theta_j}}\in H^*(N_i;\C)\quad
(R_j=\mbox{rank}_{\C}(\nu(N_i,\theta_j)))\,,
$$
where $\prod_k(1+x_k)$ is equal to the total Chern class of $\nu(N_i,\theta_j)$.  For $1\leq k\leq p-1$, 
$\varepsilon=-1,\,0,\,+1$ and an integer $\zeta$, we define numbers $T_i(k,\varepsilon,\zeta)$ and 
$S_{\varepsilon}(\zeta)$ by
\begin{eqnarray*}
&& T_i(k,\varepsilon,\zeta)
=\frac1{1-\alpha_p^k}
(\alpha_p^{k(-\varepsilon r_i+\zeta\kappa_i)}e^{-\varepsilon c_1(K_M^{-1}|_{N_i})+\zeta c_1(L|_{N_i})}-1)^{m+1}
\mbox{Td}(TN_i)\Phi(\nu(N_i,M))[N_i]\,,\\
&& S_{\varepsilon}(\zeta)=\frac1{p}\sum_{i=1}^{n}\sum_{k=1}^{p-1}T_i(k,\varepsilon,\zeta)\,,
\end{eqnarray*}
where $\mbox{Td}(TN_i)$ is the Todd class of $TN_i$.  Then $F_L(g_p)$ is defined by
\begin{eqnarray*}
&& F_L(g_p)=(m+1)\sum_{i=0}^m(-1)^i\comb{m}{i}\left(S_{-1}(m-2i)-S_{+1}(m-2i)\right)\\
&& \phantom{F_L(g_p)=}
-m\mu_{\Omega}\sum_{i=0}^{m+1}(-1)^i\comb{m+1}{i}S_{0}(m+1-2i)\,.
\end{eqnarray*}

The lifting of the character $f_{\Omega}$ given by Nakagawa is expressed by a Simons character of 
a certain foliation.  In {\rm \cite{FT}}, we gives a localization formula for the Simons character 
under Assumption \ref{asp}.  The next theorem follows from {\rm \cite[Theorem 4.7]{N2}} and 
{\rm \cite[Theorem 2.5]{FT}}.
\begin{The}\label{known} 
There exists a non-zero constant $A(m,n)$ determined only by $m,\,n$ such that
$F_L(g_p)\equiv A(m,n)f_{\Omega}(X_p)\pmod{\Z+\mu_{\Omega}\Z}$.
\end{The}
\section{Main result}
For $m,\,n\geq 1$, let $H_m,\,H_n$ be the hyperplane bundles over the complex projective spaces 
$\cp{m},\,\cp{n}$ respectively, and
$$
\pi_1:H_m\longrightarrow\cp{m}\,,\quad\pi_2:H_n\longrightarrow\cp{n}
$$
their projections.  Let $E=\pi_1^*H_m\oplus\pi_2^*H_n$ be the rank $2$ vector bundle over $\cp{m}\times\cp{n}$.  
Let $M$ be the total space of the projective bundle of $E$ and $J_M$ the tautological bundle of $M$.  
Then $M$ is an $(m+n+1)$-dimensional simply-connected compact K\"ahler manifold and the same argument as 
in {\rm \cite[Proposition 3.1]{F1}} shows that $M$ is a Fano manifold 
(see also {\rm \cite[Proposition 4.2.1]{F3}}) and the identity component of $\am{M}$ coincides with 
the factor group $(\GL{m+1}\times\GL{n+1})/\C^*$, where $\C^*$ is the center of $\GL{m+n+2}$.  
Hence the Lie algebra $\hm{M}$ is isomorphic to
$$
\{(A,B)\in\gll{m+1}\oplus\gll{n+1}\,|\,\mbox{Tr}\, A+\mbox{Tr}\, B=0\}\,,
$$
which satisfies the condition of Matsushima.

Applying the Gysin exact sequence to the fibration
$$
F=\cp{1}\longrightarrow M\stackrel{p}{\longrightarrow}B=\cp{m}\times\cp{n}\,,
$$
we have the split exact sequence
\begin{eqnarray*}
&& \hz{-1}{B}=0\,\longrightarrow\,\hz{2}{B}\simeq\hz{2}{\cp{m}}\oplus\hz{2}{\cp{n}}\\
&& \stackrel{p^*}{\longrightarrow}
\hz{2}{M}\stackrel{f}{\longrightarrow}\hz{0}{B}\simeq\Z\,\longrightarrow\,\hz{3}{B}=0\,,
\end{eqnarray*}
where $f$ is the integration along the fiber.  Then $H_m,\,H_n$ are naturally regarded as 
vector bundles over $\cp{m}\times\cp{n}$, and since $f(c_1(J_M^*))=1$, it follows that
$$
\hz{2}{M}=\{\lambda\tu+\mu\tv+\nu\tw\,|\,\lambda,\,\mu,\,\nu\in\Z\}\simeq\Z^3\,,
$$
where $\tu=c_1(p^*H_m)\,,\;\tv=c_1(p^*H_n)$ and $\tw=c_1(J_M^*)$.

\begin{rem}\label{KC}{\rm 
Let $\hat{u},\,\hat{v}$ be the first Chern forms of $H_m,\,H_n$, respectively.  Then $x\hat{u}+y\hat{v}$ 
is a K\"ahler form on $\cp{m}\times\cp{n}$ for $x,\,y>0$, and hence $x\tu+y\tv+z\tw$ is a K\"ahler class 
of $M$ for $x,\,y>0$ and sufficiently small $z>0$.  Therefore the set of K\"ahler classes of $M$ is 
contained in the subset $\{x\tu+y\tv+z\tw\,|\,x,\,y,\,z>0\}$ of $\hr{2}{M}\simeq\R^3$.}
\end{rem}

Now, let $F(x,y,z)$ be an integral homogeneous polynomial of degree $m+n+4$ defined by
$$
F(x,y,z)
=-(m(m+2)yz+n(n+2)xz+2xy)g(x,y,z)+xyz\, h(x,y,z)\,,
$$
where
\begin{eqnarray*}
&& g(x,y,z)=\sum_{s=0}^{m+n}\sum_{q=0}^{m}\comb{m+n+2}{s}\comb{s}{m-q}\comb{m+n-s}{q}(-1)^{m+n+s+q+1}\\
&& \phantom{g(x,y,z)=}
\left((x-z)^{m-q}y^{n+q+2}-x^{m-q}(y-z)^{n+q+2}\right)\,,\\
&& h(x,y,z)=\sum_{s=0}^{m+n}\sum_{q=0}^{m}\comb{m+n+2}{s}\comb{s}{m-q}\comb{m+n-s}{q}(-1)^{m+n+s+q+1}\\
&& \phantom{g(x,y,z)=}
\left(\begin{array}{l}
\{(m+n+2-s)+(n+2)(s-m+q)\}(x-z)^{m-q}y^{n+q+1}\\
\phantom{\qquad}
+m(m-q)(x-z)^{m-q-1}y^{n+q+2}\\
+\{(m+n+2-s)-n(s-m+q)\}x^{m-q}(y-z)^{n+q+1}\\
\phantom{\qquad}
-(m+2)(m-q)x^{m-q-1}(y-z)^{n+q+2}
\end{array}\right)\,.
\end{eqnarray*}
For example, if $(m,n)=(1,2)$, we have
\begin{eqnarray*}
&& F(x,y,z)=120x^2y^3z^2-420x^2y^2z^3+390x^2yz^4-120x^2z^5+60xy^4z^2-90xy^3z^3\\
&& \phantom{F(x,y,z)=}
+150xy^2z^4-99xyz^5+24xz^6-90y^4z^3+90y^3z^4-45y^2z^5+9yz^6\,.
\end{eqnarray*}

Our main result is the next theorem.
\begin{The}\label{main} 
The character $f_{\Omega}$ for $\Omega=x\tu+y\tv+z\tw$ vanishes if and only if $F(x,y,z)=0$.  
Hence the open subset of $\hr{2}{M}\simeq\R^3$ defined by $F(x,y,z)\ne 0$ does not contain any 
K\"ahler metric of constant scalar curvature. $($See Remark \ref{suff}.$)$
\end{The}
\begin{rem}{\rm 
The group $\am{M}$ contains an $(m+n+1)$-dimensional algebraic torus.  
Hence $M$ is toric and the character can be calculated also by the formula of Nakagawa {\rm \cite{N1}}.}
\end{rem}
\section{Proof of the Theorem}
Let $q\in M$, $q_m\in p^*H_m$, $q_n\in p^*H_n$ and $q_J\in J_M^*$ be points.
Then the point $q$ and the set $(q_m,q_n,q_J)$ are expressed as follows:
\begin{eqnarray*}
&& q=[(z_0,\cdots,z_m),(w_0,\cdots,w_n),(\eta_0,\eta_1)]\\
&& \phantom{q}
=[(az_0,\cdots,az_m),(bw_0,\cdots,bw_n),(ca\eta_0,cb\eta_1)]\,,\\
&& (q_m,q_n,q_J)=[[(z_0,\cdots,z_m),(w_0,\cdots,w_n),(\eta_0,\eta_1)],h_m,h_n,\xi]\\
&& \phantom{(q_m,q_n,q_J)}
=[[(az_0,\cdots,az_m),(bw_0,\cdots,bw_n),(ca\eta_0,cb\eta_1)],ah_m,bh_n,c\xi]
\end{eqnarray*}
for $a,\,b,\,c\in\C^*$.
\begin{rem}{\rm 
Since $f_{\Omega}$ vanishes on $[\hm{M},\hm{M}]$ and 
$\hm{M}/[\hm{M},\hm{M}]$ is represented by the vector field along the fiber $\cp 1$, 
the character $f_{\Omega}$ vanishes if and only if $f_{\Omega}(X)=0$ for 
the vector field $X$ along the fiber.}
\end{rem}

Now we assume that $p$ is an odd prime number hereafter.  Then an action of 
$\Z_p=\langle g_p\rangle\subset(\GL{m+1}\times\GL{n+1})/\C^*$ on $M$ is defined by
\begin{eqnarray}
&& g_p\cdot\left[(z_0,\cdots,z_m),(w_0,\cdots,w_n),(\eta_0,\eta_1)\right]
=\left[(z_0,\cdots,z_m),(\alpha_p w_0,\cdots,\alpha_p w_n),(\eta_0,\eta_1)\right]\,.
\label{gAct}
\end{eqnarray}
This action naturally extends to an inclusion $\iota:\uone\longrightarrow\am{M}$, which defines 
vector fields $X,\,X_p\in\hm{M}$ along the fiber as in (\ref{xp}) and we have $g_p=\exp(X_p)$.  
The fixed point set of $g_p^k$ has the following two connected components
$$
N_1=\left[(z_0,\cdots,z_m),(w_0,\cdots,w_n),(1,0)\right]\;,\qquad 
N_2=\left[(z_0,\cdots,z_m),(w_0,\cdots,w_n),(0,1)\right]
$$
for $1\leq k\leq p-1$, which are isomorphic to $\cp{m}\times\cp{n}$ and have cell decompositions 
with no codimension one cells.  Let $\nu(N_i,M)$ be the normal bundle of $N_i\;(i=1,\,2)$ in $M$.  
Then, since 
\begin{eqnarray*}
&& \left[(z_0,\cdots,z_m),(w_0,\cdots,w_n),(1,\tau)\right]
=\left[(az_0,\cdots,az_m),(bw_0,\cdots,bw_n),(1,a^{-1}b\tau)\right]\,,\\
&& g_p\cdot\left[(z_0,\cdots,z_m),(w_0,\cdots,w_n),(1,\tau)\right]
=\left[(z_0,\cdots,z_m),(w_0,\cdots,w_n),(1,\alpha_p^{-1}\tau)\right]\,,
\end{eqnarray*}
we have
$$
\nu(N_1,M)\simeq H_m^{-1}\otimes H_n\,,\qquad g_p|\nu(N_1,M)=g_p|(K_M^{-1}|_{N_1})=\alpha_p^{-1}\,.
$$
The same argument shows that
$$
\nu(N_2,M)\simeq H_m\otimes H_n^{-1}\,,\qquad g_p|\nu(N_2,M)=g_p|(K_M^{-1}|_{N_2})=\alpha_p\,.
$$
Hence it follows from the equality 
$c_1(K_M^{-1}|_{N_i})=c_1(M)|_{N_i}=c_1(TN_i)+c_1(\nu(N_i,M))$ that
\begin{eqnarray*}
&& c_1(\nu(N_1,M))=-u+v\;,\qquad c_1(\nu(N_2,M))=u-v\,,\\
&& c_1(K_M^{-1}|_{N_1})=mu+(n+2)v\;,\qquad 
c_1(K_M^{-1}|_{N_2})=(m+2)u+nv\,,
\end{eqnarray*}
where $u=c_1(H_m),\,v=c_1(H_n)$.  
It is obvious that $\tu|_{N_i}=u$, $\tv|_{N_i}=v$ for $i=1,\,2$.  Also, since
$$
[[(z_0,\cdots,z_m),(w_0,\cdots,w_n),(1,0)],\xi]
=[[(az_0,\cdots,az_m),(bw_0,\cdots,bw_n),(1,0)],a^{-1}\xi]\,,
$$
it follows that $\tw|_{N_1}=-u$.  The same argument shows that $\tw|_{N_2}=-v$.
Using the equalities above, we see that
$$
c_1(M)=(m+2)\tu+(n+2)\tv+2\tw\,,
$$
and hence for $\Omega=x\tu+y\tv+z\tw$ it follows that
\begin{equation}\label{mo}
\mu_{\Omega}=\frac{m(m+2)yz+n(n+2)xz+2xy}{(m+n+1)xyz}\,.
\end{equation}

Let $\lambda,\,\mu,\,\nu$ be integers.  Then $\Omega=\lambda\tu+\mu\tv+\nu\tw$ 
coincides with the first Chern class of the complex line bundle $L$ defined by
$$
L=p^*H_m^{\lambda}\otimes p^*H_n^{\mu}\otimes (J_M^*)^{\nu}\,.
$$
The action (\ref{gAct}) lifts to actions on $p^*H_m,\,p^*H_n,\,J_M^*$ as follows:
\begin{eqnarray*}
&& g_p\cdot[[(z_0,\cdots,z_m),(w_0,\cdots,w_n),(\eta_0,\eta_1)],h_m,h_n,\xi]\\
&& \phantom{g_p\cdot}
=[[(z_0,\cdots,z_m),(\alpha_p w_0,\cdots,\alpha_p w_n),(\eta_0,\eta_1)],h_m,h_n,\xi]\,.
\end{eqnarray*}
This action defines a lift of the action (\ref{gAct}) to $L$ and we can show that
\begin{eqnarray*}
&& g_p|(p^*H_m|_{N_i})=1\;,\quad g_p|(p^*H_n|_{N_i})=\alpha_p^{-1}\quad(i=1,\,2)\\
&& g_p|(J_M^*|_{N_1})=1\;,\quad g_p|(J_M^*|_{N_2})=\alpha_p\,,
\end{eqnarray*}
and hence that
\begin{equation}
g_p|(L|_{N_1})=\alpha_p^{-\mu}\;,\qquad g_p|(L|_{N_2})=\alpha_p^{-\mu+\nu}\,.
\end{equation}
Using the results above, we have
\begin{eqnarray*}
&& T_i(k,\varepsilon,\zeta)=u^mv^n\mbox{-coeff. of }\\
&& \phantom{T_i(k,\varepsilon,\zeta)=}
\frac1{1-\alpha_p^k}\left(\alpha_p^{k(-\varepsilon r+\zeta\kappa)}
e^{-\varepsilon(au+bv)+\zeta(\rho u+\tau v)}-1\right)^{m+n+2}\\
&& \phantom{T_i(k,\varepsilon,\zeta)=}
\left(\frac{u}{1-e^{-u}}\right)^{m+1}\left(\frac{v}{1-e^{-v}}\right)^{n+1}
\frac1{1-\alpha_p^{-k\delta}e^{-\delta(u-v)}}\,,
\end{eqnarray*}
where $r,\,\kappa,\,a,\,b,\,\rho,\,\tau,\,\delta$ are numbers determined by $i$ as follows:
$$
\begin{tabular}{|c|c|c|c|c|c|c|c|}
\hline
 & $r$ & $\kappa$ & $a$ & $b$ & $\rho$ & $\tau$ & $\delta$\\
\hline
$i=1$ & $-1$ & $-\mu$ & $m$ & $n+2$ & $\lambda-\nu$ & $\mu$ & $-1$\\
\hline
$i=2$ & $1$ & $-\mu+\nu$ & $m+2$ & $n$ & $\lambda$ & $\mu-\nu$ & $1$\\
\hline
\end{tabular}
$$
Then, using the substitution $x=e^{u}-1,\,y=e^{v}-1$, we have
\begin{eqnarray*}
&& T_i(k,\varepsilon,\zeta)=u^{-1}v^{-1}\mbox{-coeff. of }\\
&& \phantom{T_i(k,\varepsilon,\zeta)=}
\frac1{1-\alpha_p^k}\left(\alpha_p^{k(-\varepsilon r+\zeta\kappa)}
e^{u(\zeta\rho-\varepsilon a)}e^{v(\zeta\tau-\varepsilon b)}-1\right)^{m+n+2}\\
&& \phantom{T_i(k,\varepsilon,\zeta)=}
\left(\frac{e^{u}}{\,e^{u}-1\,}\right)^{m+1}
\left(\frac{e^{v}}{\,e^{v}-1\,}\right)^{n+1}
\frac1{1-\alpha_p^{-k\delta}e^{-\delta u}e^{\delta v}}\\
&& \phantom{T_i(k,\varepsilon,\zeta)}
=\left(\frac1{2\pi i}\right)^2\oint_{C(u)}\oint_{C(v)}
\frac1{1-\alpha_p^k}\left(\alpha_p^{k(-\varepsilon r+\zeta\kappa)}
e^{u(\zeta\rho-\varepsilon a)}e^{v(\zeta\tau-\varepsilon b)}-1\right)^{m+n+2}\\
&& \phantom{T_i(k,\varepsilon,\zeta)=\left(\frac1{2\pi i}\right)^2\oint_{C(u)}\oint_{C(v)}}
\frac{\left(e^{u}\right)^m}{\left(e^{u}-1\right)^{m+1}}
\frac{\left(e^{v}\right)^n}{\left(e^{v}-1\right)^{n+1}}
\frac1{1-\alpha_p^{-k\delta}e^{-\delta u}e^{\delta v}}e^{u}e^{v}\;dv du\\
&& \mbox{(where $C(u),\,C(v)$ are sufficiently small counterclockwise loops around the origin)}\\
&& \phantom{T_i(k,\varepsilon,\zeta)}
=\left(\frac1{2\pi i}\right)^2\oint_{C(x)}\oint_{C(y)}
\frac1{1-\alpha_p^k}\left(\alpha_p^{k(-\varepsilon r+\zeta\kappa)}
(1+x)^{\zeta\rho-\varepsilon a}(1+y)^{\zeta\tau-\varepsilon b}-1\right)^{m+n+2}\\
&& \phantom{T_i(k,\varepsilon,\zeta)=\left(\frac1{2\pi i}\right)^2\oint_{C(x)}\oint_{C(y)}}
\frac{(1+x)^m}{x^{m+1}}\frac{(1+y)^n}{y^{n+1}}
\frac1{1-\alpha_p^{-k\delta}(1+x)^{-\delta}(1+y)^{\delta}}\;dydx\\
&& \mbox{(where $C(x),\,C(y)$ are sufficiently small counterclockwise loops around the origin)}\,.
\end{eqnarray*}
Here we set $\beta=\zeta\rho-\varepsilon a\;,\quad \gamma=\zeta\tau-\varepsilon b$ and
\begin{eqnarray*}
&& \Phi=(1+x)^{-\delta}(1+y)^{\delta}-1=-\delta x+\delta y+Q(x,y)\,,\\
&& \Psi=(1+x)^{\beta}(1+y)^{\gamma}-1=\beta x+\gamma y+R(x,y)\,,
\end{eqnarray*}
where the total degrees of $Q(x,y),\,R(x,y)$ are greater than $1$.  Then we have
\begin{eqnarray*}
&& T_i(k,\varepsilon,\zeta)\\
&& =x^my^n\mbox{-coeff. of }\\
&& \phantom{=}
\frac1{\,1-\alpha_p^k\,}\left(\alpha_p^{k(\zeta\kappa-\varepsilon r)}-1+
\alpha_p^{k(\zeta\kappa-\varepsilon r)}\Psi\right)^{m+n+2}
(1+x)^m(1+y)^n\left(1-\alpha_p^{-k\delta}-\alpha_p^{-k\delta}\Phi\right)^{-1}\\
&& =x^my^n\mbox{-coeff. of }\\
&& \phantom{=}
\frac1{\,1-\alpha_p^k\,}\sum_{s=0}^{m+n}\comb{m+n+2}{s}
\left(\alpha_p^{k(\zeta\kappa-\varepsilon r)}-1\right)^{m+n+2-s}
\alpha_p^{ks(\zeta\kappa-\varepsilon r)}\Psi^s
(1+x)^m(1+y)^n\\
&& \phantom{=\frac1{\,1-\alpha_p^k\,}}
\sum_{j=0}^{m+n}\frac{\,\alpha_p^{-kj\delta}\Phi^j\,}{\,\left(1-\alpha_p^{-k\delta}\right)^{j+1}\,}\\
&& =x^my^n\mbox{-coeff. of }\\
&& \phantom{=}
\sum_{s=0}^{m+n}\sum_{j=0}^{m+n-s}
\comb{m+n+2}{s}(-1)\Lambda_j(\alpha_p^k)(1+x)^m(1+y)^n\Phi^j\Psi^s\,,
\end{eqnarray*}
where $\Lambda_j(t)$ is an element of $\Z[t,t^{-1}]$ defined by
$$
\Lambda_j(t)=\frac{\,t^{s(\zeta\kappa-\varepsilon r)+\delta}
\left(t^{\zeta\kappa-\varepsilon r}-1\right)^{m+n+2-s}\,}{\,(t-1)(t^{\delta}-1)^{j+1}\,}\,.
$$
Here, since
$$
\sum_{k=1}^{p-1}\alpha_p^{kl}\equiv -1\pmod{p}
$$
for any integer $l$, we have
\begin{eqnarray*}
&& (-1)\sum_{k=1}^{p-1}\Lambda_j(\alpha_p^k)\equiv \Lambda_j(1)\pmod{p}\\
&& \phantom{(-1)\sum_{k=1}^{p-1}\Lambda_j(\alpha_p^k)}
=\left\{\begin{array}{cl} 
0 & \mbox{if }j<m+n-s\\
\delta^{m+n-s+1}(\zeta\kappa-\varepsilon r)^{m+n+2-s} & \mbox{if }j=m+n-s
\end{array}\right.\,.
\end{eqnarray*}
Therefore we have
\begin{eqnarray*}
&& \sum_{k=1}^{p-1}T_i(k,\varepsilon,\zeta)\\
&& \equiv\;x^my^n\mbox{-coeff. of }\\
&& \phantom{\equiv\;}
\sum_{s=0}^{m+n}\comb{m+n+2}{s}\delta^{m+n-s+1}(\zeta\kappa-\varepsilon r)^{m+n+2-s}
(-\delta(x-y))^{m+n-s}(\beta x+\gamma y)^s\pmod{p}\\
&& =\;x^my^n\mbox{-coeff. of }\\
&& \phantom{\equiv\;}
\sum_{s=0}^{m+n}\comb{m+n+2}{s}\delta^{m+n-s+1}(\zeta\kappa-\varepsilon r)^{m+n+2-s}(-\delta)^{m+n-s}\\
&& \phantom{\equiv\;}
\sum_{h=0}^s\comb{s}{h}\beta^h x^h\gamma^{s-h}y^{s-h}\sum_{q=0}^{m}\comb{m+n-s}{q}x^q(-y)^{m+n-s-q}\\
&& =\sum_{s=0}^{m+n}\sum_{q=0}^m\comb{m+n+2}{s}\comb{s}{m-q}\comb{m+n-s}{q}(-1)^q\\
&& \phantom{\sum_{k=1}^{p-1}S_1(k)=}
\delta(\kappa\zeta-r\varepsilon)^{m+n+2-s}(\rho\zeta-a\varepsilon)^{m-q}(\tau\zeta-b\varepsilon)^{s-m+q}\,,
\end{eqnarray*}
and hence it follows that
\begin{eqnarray*}
&& S_{\varepsilon}(\zeta)\\
&& \equiv\frac1p\sum_{s=0}^{m+n}\sum_{q=0}^m\comb{m+n+2}{s}\comb{s}{m-q}\comb{m+n-s}{q}(-1)^q\\
&& \phantom{=\frac1p}
\left(\begin{array}{l} 
(-1)^{m+n+s+1}(\mu\zeta-\varepsilon)^{m+n+2-s}((\lambda-\nu)\zeta-m\varepsilon)^{m-q}
(\mu\zeta-(n+2)\varepsilon)^{s-m+q}\\
+((-\mu+\nu)\zeta-\varepsilon)^{m+n+2-s}(\lambda\zeta-(m+2)\varepsilon)^{m-q}((\mu-\nu)\zeta-n\varepsilon)^{s-m+q}
\end{array}\right)\pmod{\Z}\\
&& =\frac1p g(\lambda,\mu,\nu)\zeta^{m+n+2}-\varepsilon \frac1p h(\lambda,\mu,\nu)\zeta^{m+n+1}+\varphi(\zeta)\,,
\end{eqnarray*}
where the degree of $\varphi(\zeta)$ is less than $m+n+1$.

Here for $f(x)=(\sinh x)^k$ we have
$$
f(x)=\frac1{\,2^k\,}\sum_{i=0}^{k}(-1)^i\comb{k}{i}e^{(k-2i)x}\;,\quad 
f(x)=x^k+\frac{k}{6}x^{k+2}+\mbox{higher order terms}
$$
and hence it follows that
$$
2^kf^{(l)}(0)
=\sum_{i=0}^{k}(-1)^i\comb{k}{i}(k-2i)^{l}
=\left\{\begin{array}{cl} 
0 & \mbox{if }0\leq l<k\mbox{ or }l=k+1\\
2^kk! & \mbox{if }l=k
\end{array}\right.\,.
$$
Therefore it follows from (\ref{mo}) that
\begin{eqnarray*}
&& \lambda\mu\nu F_L(g_p)\\
&& =(m+n+2)\lambda\mu\nu\\
&& \phantom{=}
\sum_{i=0}^{m+n+1}(-1)^i\comb{m+n+1}{i}\left(S_{-1}(m+n+1-2i)-S_{+1}(m+n+1-2i)\right)\\
&& \phantom{=}
-(m(m+2)\mu\nu+n(n+2)\lambda\nu+2\lambda\mu)\sum_{i=0}^{m+n+2}(-1)^i\comb{m+n+2}{i}S_{0}(m+n+2-2i)\\
&& \equiv\frac{2^{m+n+2}(m+n+2)!}{p} F(\lambda,\mu,\nu)\pmod{\Z}\,.
\end{eqnarray*}
Hence, for any odd prime number $p$, it follows from Theorem \ref{known} that
\begin{eqnarray*}
&& \frac1p A(m,n)\lambda\mu\nu f_{\Omega(\lambda,\mu,\nu)}(X)
=A(m,n)\lambda\mu\nu f_{\Omega(\lambda,\mu,\nu)}(X_p)\\
&& \phantom{\frac1p A(m,n)\lambda\mu\nu f_{\Omega(\lambda,\mu,\nu)}(X)}
\equiv\frac1p 2^{m+n+2}(m+n+2)!F(\lambda,\mu,\nu)\pmod{\Z}\,,
\end{eqnarray*}
where $\Omega(\lambda,\mu,\nu)=\lambda\tu+\mu\tv+\nu\tw$, which implies that
\begin{equation}\label{equal}
A(m,n)\lambda\mu\nu f_{\Omega(\lambda,\mu,\nu)}(X)=2^{m+n+2}(m+n+2)!F(\lambda,\mu,\nu)\,.
\end{equation}

Now, since $\triangle_{k\omega}=k^{-1}\triangle_{\omega}$, it follows that 
$xyz f_{\Omega(x,y,z)}(X)$ is a homogeneous function in $x,\,y,\,z$ of degree $m+n+4$ 
as well as $F(x,y,z)$.  Moreover, since the set
$$
\{(x,y,z)\in\R^3\,|\,(rx,ry,rz)\in\Z^3\;\mbox{for some}\;r>0\}
$$
is dense in $\R^3$, the equality (\ref{equal}) implies that for any $(x,y,z)\in\R^3$
$$
A(m,n)xyz f_{\Omega(x,y,z)}(X)=2^{m+n+2}(m+n+2)!F(x,y,z)\,.
$$
The result in Theorem \ref{main} follows immediately from the equality above.
\begin{rem}\label{suff}{\rm 
Let $G=(U(m+1)\times U(n+1))/U(1)$ be the maximal compact subgroup of $\am{M}$ and 
$q=[(z_0,\cdots,z_m),(w_0,\cdots,w_n),(\eta_0,\eta_1)]$ a point in $M$.  
Then we can see that the real dimension of the isotropy subgroup of $G$ at $q$ is equal to 
$m^2+n^2$ if $\eta_0\eta_1\ne 0$ and is equal to $m^2+n^2+1$ if $\eta_0\eta_1=0$, 
which implies that the real codimension of the principal orbit of $G$ in $M$ is one.  
Hence it follows from Corollary 1.1 in {\rm \cite{H}} that each K\"ahler class of $M$ contains an extremal 
metric, and therefore it follows from {\rm \cite[Theorem 4]{C}} (see also {\rm \cite[Theorem 3.3.1]{F3}}) 
that a K\"ahler class contains a K\"ahler metric of constant scalar curvature if the character for the 
K\"ahler class vanishes.  Hence a K\"ahler class $\Omega=x\tu+y\tv+z\tw$ contains a K\"ahler metric of 
constant scalar curvature if and only if $F(x,y,z)=0$.  Moreover we can see that the $\am{M}$-orbit of 
$q$ with $\eta_0\eta_1\ne 0$ coincides with the open subset $M\setminus(N_1\cup N_2)$ of $M$.  
Hence $M$ is an almost-homogeneous manifold (see {\rm \cite{HS}}) and therefore it follows from 
{\rm \cite[Theorem 4]{H}} that $M$ admits a K\"ahler metric of constant scalar curvature.}
\end{rem}
\section{Examples}
In this section, we consider the cases $1\leq m<n\leq 10$.  
Since $F(x,y,z)$ is a homogeneous polynomial, $F(x,y,z)$ for $x,\,y,\,z>0$ is determined by its restriction 
to the face $f$ of a regular octahedron defined by
$$
f=\{(x,y,z)\,|\,x+y+z=1\,,\;x,\,y,\,z>0\}\,.
$$
Let $C$ be a point in $f$ defined by
$$
C=\frac1{m+n+6}(m+2,n+2,2)
$$
and set $A=(1,0,0),\,B=(0,1,0)$.  Then, since $C$ is homothetic to $c_1(M)>0$, $C$ is a K\"ahler class 
and hence the interior of the triangle ABC is contained in the set of K\"ahler classes of $M$ 
(see Remark \ref{KC}).  Let $l_1,\,l_2$ be lines in $f$ defined by
\begin{eqnarray*}
&& l_1(t)=(x_1(t),y_1(t),z_1(t))=(1-t)A+tC\;,\\
&& l_2(t)=(x_2(t),y_2(t),z_2(t))=(1-t)(\frac12,\frac12,0)+t(0,0,1)
\end{eqnarray*}
for $0<t<1$.  Then we have
\begin{eqnarray*}
&& \lim_{t\to+0}F(l_1(t))/y_1(t)^{n+3}\\
&& =\lim_{t\to+0}\sum_{s=0}^{m+n}\sum_{q=0}^{m}\comb{m+n+2}{s}\comb{s}{m-q}\comb{m+n-s}{q}
(-1)^{m+n+s+q+1}t^q\\
&& \phantom{\kern 2.8cm}
2(n+2)^{-n-2}(m+n+6)^{-q}\\
&& \phantom{\kern 2.8cm}
\left\{((n+1)s+(n+2)q-m-mn-2n-n^2)(n+2)^{n+1+q}\right.\\
&& \phantom{\kern 4cm}
\left.
-((n+1)s+nq-m-mn-2n-n^2-2)n^{n+1+q}\right\}\\
&& =\sum_{s=0}^{m+n}\comb{m+n+2}{s}\comb{s}{m}(-1)^{m+n+s+1}2(n+2)^{-n-2}\\
&& \phantom{\kern 1.2cm}
\left\{((n+1)s-m-mn-2n-n^2)(n+2)^{n+1}\right.\\
&& \phantom{\kern 1.4cm}
\left.
-((n+1)s-m-mn-2n-n^2-2)n^{n+1}\right\}\,,\\
&& \lim_{t\to+0}F(l_2(t))/z_2(t)^2\\
&& =\sum_{s=0}^{m+n}\sum_{q=0}^{m}\comb{m+n+2}{s}\comb{s}{m-q}\comb{m+n-s}{q}
(-1)^{m+n+s+q+1}2^{-m-n-2}\\
&& \phantom{=\sum_{s=0}^{m+n}\sum_{q=0}^{m}}
\left\{2(n-m)q^2+\left(2(n+1)s+(m^2-4mn-n^2-7m-3n-2)\right)q\right.\\
&& \phantom{=\sum_{s=0}^{m+n}\sum_{q=0}^{m}\left\{\right.}
\left.+(-mn+n^2-m+2n+1)s+3m^2+m^2n-n^3-mn-4n^2-2m-4n\right\}\,.
\end{eqnarray*}
Direct computation using the equalities above shows that
$$
\lim_{t\to+0}F(l_1(t))/y_1(t)^{n+3}<0\;,\qquad 
\lim_{t\to+0}F(l_2(t))/z_2(t)^2>0\,,
$$
which imply that there exist points $P_1,\,P_2$ in the interior of the triangle ABC 
such that $F(P_1)<0$ and $F(P_2)>0$.  
Therefore there exist infinitely many K\"ahler classes $\Omega$ such that $f_{\Omega}$ vanishes 
and hence that $\Omega$ contains a K\"ahler metric of constant scalar curvature 
(see Remark \ref{suff}).

On the other hand, direct computation also shows that
$$
F(m+2,n+2,2)\ne 0\,,
$$
which implies that $c_1(M)$ does not contain any K\"ahler metric of constant scalar curvature.  
This result shows that $M$ does not admit any K\"ahler-Einstein metric.  (See {\rm \cite{F1}}.)

\quad\\
\begin{flushleft}
\textsc{
Tokyo University of Marine science and technology \\
4-5-7 Kounan, Minato-ku \\
Tokyo 108-8477 Japan}\\
\textit{E-mail address}: 
tsubois@kaiyodai.ac.jp
\end{flushleft}

\begin{thebibliography}{99}
%
%
%
\bibitem{B} S. Bando, An obstruction for Chern class forms to be harmonic, Kodai Math. J. {\bf 29} 
(2006), 337--345.

\bibitem{C} E. Calabi, Extremal K\"ahler metrics II, Differential geometry and complex analysis, 
(I. Chavel and H.M. Farkas eds.), 95--114, Springer-Verlag, Berline-Heidelberg-New York, 1985.

\bibitem{F1} A. Futaki, An obstruction to the existence of Einstein-K\"ahler metrics, Invent. Math. {\bf 73} 
(1983), 437--443.

\bibitem{F2} A. Futaki, On compact K\"ahler manifold of constant scalar curvature, Proc. Japan Acad. Ser. A 
{\bf 59} (1983), 401--402.

\bibitem{F3} A. Futaki, K\"ahler-Einstein metrics and integral invariants, 
Lecture Notes in Math. {\bf 1314}, Springer-Verlag, Berlin, 1988.

\bibitem{FM} A. Futaki and S. Morita, Invariant polynomials of the automorphism group of a compact complex manifold, 
J. Differential Geom. {\bf 21}(1985), 135--142.

\bibitem{FT} A. Futaki and K. Tsuboi, Fixed point formula for characters of automorphism groups associated with 
K\"ahler classes, Math. Res. Lett. {\bf 8} (2001), 495--507.

\bibitem{H} A. D. Hwang, On existence of K\"ahler metrics with constant scalar curvature, 
Osaka J. Math. {\bf 31} (1994), 561--595.

\bibitem{HS} A. T. Huckleberry and D. M. Snow, Almost-homogeneous K\"ahler manifolds with hypersurface orbits, 
Osaka J. Math. {\bf 19} (1982), 763--786.

\bibitem{KW} J. Kazdan and F. Warner, Prescribing curvatures, Proc. Sympos. Pure Math. {\bf 27} (1975), 309--319.

\bibitem{K} S. Kobayashi, Transformation groups in differential geometry, 
Springer-Verlag, Berlin-Heidelberg-New York, 1972.

\bibitem{L1} A. Lichnerowicz, Sur les transformations analytiques d'une vari\'et\'e 
K\"ahlerienne compacte, 1959, Colloque Geom. Diff. Global (Bruxelles, 1958), 11--26, 
Centre Belge Rech. Math., Louvain.

\bibitem{L2} A. Lichnerowicz, Isom\'etrie et transformations analytiques d'une vari\'et\'e 
K\"ahlerienne compacte, Bull. Soc. Math. France {\bf 87} (1959), 427--437.

\bibitem{M} Y. Matsushima, Sur la structure du groupe d'hom\'eomorphismes d'une certaine 
vari\'et\'e Kaehl\'erienne, Nagoya Math. J. {\bf 11} (1957), 145--150.

\bibitem{N1} Y. Nakagawa, Bando-Calabi-Futaki character of compact toric manifolds, 
Tohoku Math. J. {\bf 53} (2001), 479--490.

\bibitem{N2} Y. Nakagawa, The Bando-Calabi-Futaki character and its lifting to a group character, 
Math. Ann. {\bf 325} (2003), 31--53.

\bibitem{T} G. Tian, K\"ahler-Einstein metrics on algebraic manifolds, in: Proc. C.I.M.E. conference on 
Transcendental methods in algebraic geometry, Lecture Notes in Math. 1646, Springer-Verlag, 
Berlin-Heidelberg-New York, 1996, 143--185.
%
\end{thebibliography}
\end{document}